%&amstex          
\input amstex\documentstyle{amsppt}  
\pagewidth{12.5cm}\pageheight{19cm}\magnification\magstep1
\topmatter
\title Unipotent elements in small characteristic, IV\endtitle
\author G. Lusztig\endauthor
\address{Department of Mathematics, M.I.T., Cambridge, MA 02139}\endaddress
\thanks{Supported in part by the National Science Foundation}\endthanks
\endtopmatter   
\document
\define\Ann{\text{\rm Ann}}

\define\be{\bar e}

\define\frl{\forall}
\define\pe{\perp}
\define\si{\sim}

\define\sqc{\sqcup}

\define\qua{\quad}

\define\dx{\dot x}

\define\bA{\bar A}
\define\baf{\bar f}

\define\op{\oplus}

\define\part{\partial}

\define\n{\notin}

\define\m{\mapsto}
\define\do{\dots}

\define\sub{\subset}    

\define\T{\times}
\define\ti{\tilde}
\define\nl{\newline}
\redefine\i{^{-1}}

\define\Ad{\text{\rm Ad}}

\define\End{\text{\rm End}}

\define\tr{\text{\rm tr}}

\define\a{\alpha}
\redefine\b{\beta}

\redefine\d{\delta}
\define\e{\epsilon}

\define\io{\iota}

\define\s{\sigma}

\define\z{\zeta}
\define\x{\xi}

\redefine\D{\Delta}

\define\Si{\Sigma}

\define\Ps{\Psi}

\define\kk{\bold k}

\define\CC{\bold C}

\define\FF{\bold F}

\define\NN{\bold N}

\define\ZZ{\bold Z}

\define\cn{\Cal N}
\define\co{\Cal O}

\define\cu{\Cal U}

\define\cy{\Cal Y}

\define\fb{\frak b}

\define\fg{\frak g}

\define\fD{\frak D}

\define\fF{\frak F}

\define\fQ{\frak Q}

\define\fS{\frak S}

\define\fU{\frak U}

\define\tm{\ti m}

\define\tV{\ti V}

\define\bQ{\bar Q}

\define\KW{KW}
\define\LI{L1}
\define\LII{L2}
\define\LIII{L3}
\define\SPR{S}
\define\XU{X}

\head Introduction\endhead
Let $\kk$ be an algebraically closed field of characteristic exponent $p\ge1$. Let $G$ be a connected reductive 
algebraic group over $\kk$ and let $\fg$ be the Lie algebra of $G$. Note that $G$ acts on $G$ and on $\fg$ by the
adjoint action and on $\fg^*$ by the coadjoint action. (For any $\kk$-vector space $V$ we denote by $V^*$ the dual
vector space.) Let $G_\CC$ be the reductive group over $\CC$ of the same type as $G$. 
Let $\cu_G$ be the variety of unipotent elements of $G$. Let $\cn_\fg$ be the variety of nilpotent elements of 
$\fg$. Let $\cn_{\fg^*}$ be the variety of nilpotent elements of $\fg^*$ (following \cite{\KW} we say that a
linear form $\x:\fg@>>>\kk$ is nilpotent if its kernel contains some Borel subalgebra of $\fg$). In 
\cite{\LI, \LII, \LIII} we have proposed a definition of a partition of $\cu_G$ and of $\cn_{\fg}$ into smooth 
locally closed $G$-stable pieces which are indexed by the unipotent classes in
$G_\CC$ and which in many ways depend very smoothly on $p$. In this paper we 
propose a definition of an analogous partition of $\cn_{\fg^*}$ into pieces which are indexed by the unipotent classes in $G_\CC$. 
(This definition is only of interest for $p>1$, small; for $p=1$ or $p$
large we can identify $\cn_{\fg}$ with $\cn_{\fg^*}$ and the partition of $\cn_{\fg^*}$ is deduced from the
partition of $\cn_{\fg}$.) We will illustrate this in the case where $G$ is of type $A,C$ or $D$ and $p$ is arbitrary.

{\it Notation.} If $f$ is a permutation of a set $X$ we denote $X^f=\{x\in X;f(x)=x\}$. The cardinal of a finite 
set $X$ is denoted by $|X|$. For any subspace $U$ of $\fg$ let $\Ann(U)=\{\x\in\fg^*;\x|_U=0\}$.

\head 1.\endhead
\subhead 1.1\endsubhead
Let $V$ be a $\kk$-vector space of finite dimension. Let $G=GL(V)$. We have $\fg=\End(V)$. We have an isomorphism

(a) $\fg@>\si>>\fg^*$
\nl
given by $X\m[T\m\tr(TX,V)]$.

\subhead 1.2\endsubhead
Let $V$ be a $\kk$-vector space of finite even dimension with a fixed nondegenerate quadratic form $Q:V@>>>\kk$.
Let $(,)$ be the (nondegenerate) symmetric bilinear form $V\T V@>>>\kk$ given by $(x,y)=Q(x+y)-Q(x)-Q(y)$ for 
$x,y\in V$. Let $G=SO(V)$ be the special orthogonal group of $Q$. We have
$\fg=\{T\in\End(V);(Tx,x)=0\qua \frl x\in V\}$. Let $\fS(V)$ be the vector space consisting of all symplectic 
forms $V\T V@>>>\kk$. The following result is easily verified.

(a) {\it We have a vector space isomorphism $\fg@>\si>>\fS(V)$, $T\m[x,y\m(Tx,y)]$. This is compatible with the 
$SO(V)$-actions where $SO(V)$ acts on $\fS(V)$ by the restriction of the obvious $GL(V)$-action.}
\nl
By taking transpose we obtain an isomorphism $\fS(V)^*@>\si>>\fg^*$ compatible with the natural $SO(V)$-actions. 
We can find an isomorphism $\fS(V^*)@>\si>>\fS(V)^*$ compatible with the natural $GL(V)$-actions and, by 
restriction, with the natural $SO(V)$-actions. Now $(,)$ defines an isomorphism $V@>>>V^*$ hence an isomorphism 
$\fS(V)@>\si>>\fS(V^*)$ compatible with the natural $SO(V)$-actions. The composition 
$\fg@>\si>>\fS(V)@>\si>>\fS(V^*)@>\si>>\fS(V)^*@>\si>>\fg^*$ is an isomorphism

(b) $\fg@>\si>>\fg^*$
\nl
compatible with the $SO(V)$-actions.

\subhead 1.3\endsubhead
Let $V$ be a $\kk$-vector space of finite (even) dimension with a fixed nondegenerate symplectic form 
$(,):V\T V@>>>\kk$. Let $G=Sp(V)$ be the symplectic group of $(,)$. We have 
$\fg=\{T\in\End(V);(Tx,y)+(x,Ty)=0\qua \frl x,y\in V\}$. Let $\fQ(V)$ be the vector space consisting of all 
quadratic forms $V@>>>\kk$. According to T. Xue \cite{\XU}:

(a) {\it we have a natural vector space isomorphism $\s_V:\fg^*@>\si>>\fQ(V)$.}
\nl
Indeed, let $Z=\{X\in\End(V);(Xa,a)=0\frl a\in V\}$. We have a diagram $\fQ(V)@<\a<<\End(V)/Z@>\b>>\fg^*$ ($\a$ is
induced by $X\m[a\m(Xa,a)]$; $\b$ is induced by $X\m[T\m\tr(TX,V)]$). Now $\a,\b$ are isomorphisms and we set
$\s_V=\a\b\i$.

We define a linear map $\fQ(V)@>>>\fg,Q\m A_Q$ by 

$Q(x+y)-Q(x)-Q(y)=(A_Qx,y)$ for all $x,y\in V$. 
\nl
If $p=2$, then for $Q\in\fQ(V)$ we have $A_Q\in\fQ'(V)$ where $\fQ'(V)=\{A\in\End(V);(Ax,x)=0\qua\frl x\in V\}$.
Let $\fQ(V)_{nil}=\{Q\in\fQ(V);A_Q\text{ is nilpotent}\}$. 

\subhead 1.4\endsubhead
In this subsection we assume that $\kk$ is an algebraic closure of a finite field $\FF_q$ and that we are given an
$\FF_q$-rational structure on $G$. Then $\fg,\fg^*,\cn_{\fg},\cn_{\fg^*}$ have induced $\FF_q$-structures.
For any $\kk$-variety $X$ with an $\FF_q$-structure we denote by $F:X@>>>X$ the corresponding Frobenius map. Let 
$N$ be the number of roots of $G$. According to \cite{\SPR} we have

(a) $|\cn_{\fg}^F|=q^N$.
\nl
We now state the following result.

(b) {\it If the adjoint group of $G$ is simple of type $A,C$ or $D$, then $|\cn_{\fg^*}^F|=q^N$.}
\nl
We can assume that $G$ is a general linear group, a symplectic group or an even special orthogonal group. The 
proof in these cases will be given in the remainder of this section. We will show elsewhere that (b) holds 
without any assumption on the type of $G$.

\subhead 1.5\endsubhead
We preserve the setup of 1.4. Assume that there exists a $G$-equivariant vector space isomorphism 
$\io:\fg@>\si>>\fg^*$ compatible with the $\FF_q$-structures. It is easy to see that $\io$ restricts to a 
bijection $\cn_\fg@>\si>>\cn_{\fg^*}$. Hence $|\cn_\fg^F|=|\cn_{\fg^*}^F|$ and 1.4(b) follows from 1.4(a). In 
particular, if $G$ is a general linear group or an even special orthogonal group or a symplectic group (with 
$p\ne2$) then $\io$ as above exists (see 1.1(a), 1.2(b)) and 1.4(b) holds in these cases.

\subhead 1.6\endsubhead
In this subsection we assume that $V,(,),G,\fg$ are as in 1.3. We set $2r=\dim V$.

Assume first that $p=2$. Let $\x\in\fg^*$, let $Q\in\fQ(V)$ be the element corresponding to $\x$ under 1.3(a) and
let $A=A_Q$. The following result is due to T. Xue \cite{\XU}.

(a) {\it If $\x\in\cn_{\fg^*}$ then $A:V@>>>V$ is nilpotent.}
\nl
Let $H_r=\{i\in\ZZ;-2r+1\le i\le 2r-1],i=\text{odd}\}$. A basis $(e_i)_{i\in H_r}$ of $V$ is said to be good if 
$(e_i,e_j)=\d_{i+j}=0$ for all $i,j\in H_r$. We can find a Borel subalgebra $\fb$ of $\fg$ such that 
$\x|_{\fb}=0$. We can find a good basis $(e_i)_{i\in H_r}$ of $V$ such that $\fb$ consists of all $T\in\End(V)$
with $Te_i=\sum_{j\in H_r;i\le j}t_{ij}e_j$ for all $i\in H_r$ and $t_{ij}+t_{-j,-i}=0$ for all $i,j\in H_r$  
$(t_{ij}\in\kk)$. We can find $X\in\End(V)$ such that $\tr(XT,V)=\x(T)$ for all $T\in\fg$. Define $x_{ij}\in\kk$
by $Xe_i=\sum_{j\in H_r}x_{ij}e_j$ for all $i\in H_r$. Then $\sum_{i,j\in H_r}x_{ij}t_{ji}=0$ for any 
$t_{ij}$, ($j\le i$) such that $t_{ij}+t_{-j,-i}=0$ for all $i,j$. It follows that $x_{ij}+x_{-j,-i}=0$ for any 
$j\le i$ and $x_{i,-i}=0$ for any $i\ge0$. Define $y_{ij}\in\kk$ by $Ae_i=\sum_{j\in H_r}y_{ij}e_j$ for all 
$i\in H_r$. We have $(Ae_i,e_j)=(Xe_i,e_j)+(e_i,Xe_j)$, $Q(e_i)=(Xe_i,e_i)$. Hence $y_{ij}=x_{ij}+x_{-j,-i}$,
$Q(e_i)=x_{i,-i}$. For $j\le i$ we have $x_{ij}+x_{-j,-i}=0$ hence $y_{ij}=0$. Thus $A$ is nilpotent and (a) is 
proved.

We show a converse to (a):

(b) {\it If $A:V@>>>V$ is nilpotent then $\x\in\cn_{\fg^*}$.}
\nl
It is enough to verify the following statement:

(c) {\it Let $Q\in\fQ(V)$. Let $A=A_Q$. Assume that $A:V@>>>V$ is nilpotent. Then there exists a good basis 
$(e_i)_{i\in H_r}$ of $V$ such that $Ae_i=\sum_{j\in H_r;i<j}y_{ij}e_j$ for all $i\in H_r$ ($y_{ij}\in\kk$) and 
$Q(e_i)=0$ for $i\ge0$.}
\nl
(Indeed if (c) holds then as in the proof of (a) we can define a Borel subalgebra $\fb$ in terms of $(e_i)$ and we
have $\x|_{\fb}=0$.)

We prove (c) by induction on $r$. When $r=0$ the result is trivial. Now assume that $r\ge1$. Since $x,y\m(Av,x)$ 
is a symplectic form on an even dimensional vector space, its radical has even dimension. Thus $\dim\ker A$ is 
even. Since $A$ is nilpotent, its kernel is $\ne0$ hence it has dimension $\ge2$. Now a quadratic form on a 
$\kk$-vector space of dimension $\ge2$ vanishes at some non-zero vector. Thus there exists $v\in V-\{0\}$ such 
that $Av=0$, $Q(v)=0$. Let $(\kk v)^\pe=\{v'\in V;(v',v)=0\}$. Let $V'=(\kk v)^\pe/\kk v$. Then $V'$ inherits a 
nondegenerate symplectic form $(,)'$ from $(,)$, a quadratic form $Q'$ from $Q$ and a nilpotent endomorphism $A'$
from $A$. Note that $A'=A_{Q'}$. By the induction hypothesis there exists a good basis $(e'_i)_{i\in H_{r-1}}$ 
(relative to $(,)'$) such that $A'e'_i=\sum_{j\in H_{r-1};i<j}y'_{ij}e'_j$ for all $i\in H_{r-1}$ 
($y'_{ij}\in\kk$) and $Q'(e'_i)=0$ for $i\ge0$. For $i\in H_{r-1}$ we denote by $e_i$ a representative of $e'_i$ 
in $(\kk v)^\pe$. We have $Ae_i=\sum_{j\in H_{r-1};i<j}y'_{ij}e_j+c_iv$ for all $i\in H_{r-1}$  ($c_i\in\kk$) and
$Q(e_i)=0$ for all $i\in H_{r-1}$, $i\ge0$. We set $e_{2r-1}=v$. Let $e_{-2r+1}$ be the unique vector in $V$ such
that $(e_{-2r+1},e_j)=\d_{2r-1,j}$. We have $Ae_{-2r+1}=\sum_iy_{-2r+1,i}e_i$ with $y_{-2r+1,i}\in\kk$. Since 
$\tr(A,V)=0$ we have $y_{-2r+1,-2r+1}=0$. Thus $(e_i)_{i\in H_r}$ has the required properties.

From (a),(b) we see that $\s_V:\fg^*@>\si>>\fQ(V)$ (see 1.3(a)) restricts to a bijection 

(d) $\cn_{\fg^*}@>\si>>\fQ(V)_{nil}$.
\nl
Note that (d) holds also when $p\ne2$ (with a simpler proof).

\subhead 1.7\endsubhead
We preserve the setup of 1.6 with $p=2$. 
We assume that $\kk,\FF_q$ are as in 1.4 and that we are given an $\FF_q$-structure on $V$ compatible with $(,)$.
Then $G,\fg,\fg^*,\fQ(V),\fQ'(V)$ have natural $\FF_q$-structures with Frobenius 
maps denoted by $F$. Note that $\fQ(V)@>>>\fQ'(V),Q\m A_Q$ induces a map $\fQ(V)^F@>>>\fQ'(V)^F$ with fibres of 
cardinal $q^{2r}$. From 1.6 we see that 
$$|\cn_{\fg^*}^F|=|\fQ(V)_{nil}^F|=q^{2r}|\{A\in\fQ'(V)^F;A\text{ nilpotent}\}|=q^{2r}q^{2r^2-2r}=q^{2r^2}.$$
(The third equality follows from 1.4(a) applied to an even special orthogonal group.) This proves 1.4(b) in our 
case.

\head 2.\endhead
\subhead 2.1\endsubhead
Let $\d\in\fD_G$ (see \cite{\LIII, 1.1}). Let $\fg=\op_{i\in\ZZ}\fg_i^\d$ be the corresponding grading of $\fg$
(see \cite{\LIII, 1.2}). For $j\in\ZZ$ let $\fg^{*\d}_j=\Ann(\op_{i;i\ne-j}\fg^\d_i)$. We have 
$\fg^*=\op_{j\in\ZZ}\fg^{*\d}_j$. Let $G^\d_{\ge0}$ be as in \cite{\LIII, 1.2}. 

Since $\fg^\d_{\ge-j+1}$ is $G^\d_{\ge0}$-stable we see that $\fg^{*\d}_{\ge j}$ is $G^\d_{\ge0}$-stable. For any
$\x\in\fg^*$ let $G_\x$ be the stabilizer of $\x$ in $G$ for the coadjoint action. Let
$$\fg_2^{*\d!}=\{\x\in\fg_2^{*\d};G_\x\sub G^\d_{\ge0}\}.$$

Let $\D\in D_G$ (see \cite{\LIII, 2.1}). As in \cite{\LIII, 2.1} we write $G_{\ge0}^\D$, $\fg^\D_{\ge i}$ 
($i\in\NN$) instead of $G^\d_{\ge0},\fg^\d_{\ge i}$ (see \cite{\LIII, 1.2}) where $\d\in\D$. For $j\in\NN$ we set
$$\fg^{*\D}_{\ge j}=\Ann(\fg^\D_{\ge-j+1}).$$
We have also $\fg^{*\D}_{\ge j}=\op_{j'\in\ZZ;j'\ge j}\fg^{*\d}_{j'}$. Since $\fg^\D_{\ge-j+1}$ is 
$G^\D_{\ge0}$-stable, we see that $\fg^{*\D}_{\ge j}$ is $G^\D_{\ge0}$-stable.

For any $\d\in\D$ we have an obvious isomorphism $\fg^{*\d}_2@>\si>>\fg^{*\D}_{\ge2}/\fg^{*\D}_{\ge3}$. Via this 
isomorphism the subset $\fg_2^{*\d!}$ of $\fg_2^{*\d}$ can be viewed as a subset $\Si^{*\d}$ of 
$\fg_{\ge2}^{*\D}/\fg_{\ge3}^{*\D}$. As in \cite{\LIII, 2.3} we see that $\Si^{*\d}$ is independent of the choice 
of $\d$ in $\D$; we will denote it by $\Si^{*\D}$. Note that $\Si^{*\D}$ is a subset of 
$\fg_{\ge2}^{*\D}/\fg_{\ge3}^{*\D}$ stable under the action of $G^\D_{\ge0}$. 

Let $\s^{*\D}\sub\fg_{\ge2}^{*\D}$ be the inverse image of $\Si^{*\D}$ under the obvious map 
$\fg_{\ge2}^{*\D}@>>>\fg_{\ge2}^{*\D}/\fg_{\ge3}^{*\D}$. Now $\s^{*\D}$ is stable under the coadjoint action of 
$G^\D_{\ge0}$ on $\fg_{\ge2}^{*\D}$ and $\x\m\x$ is a map
$$\Ps_{\fg^*}:\sqc_{\D\in D_G}\s^{*\D}@>>>\cn_{\fg^*}.$$

\proclaim{Theorem 2.2} Assume that the adjoint group of $G$ is a product of simple groups of type $A,C,D$. Then 
$\Ps_{\fg^*}$ is a bijection.
\endproclaim
The general case reduces easily to the case where $G$ is almost simple of type $A,C$ or $D$. Moreover we can 
assume that $G$ is a general linear group, a symplectic group or an even special orthogonal group. The proof in 
these cases will be given in 2.4, 2.13. We expect that the theorem holds without restriction on $G$.

\subhead 2.3\endsubhead
In this subsection we assume that there exists a $G$-equivariant vector space isomorphism $\io:\fg@>\si>>\fg^*$.
Let $\d\in\fD_G$. For any $i\in\ZZ$ we have 

(a) $\io(\fg_i^\d)=\fg_i^{*\d}$. 
\nl
Recall that $\fg_i^\d=\{x\in\fg;\Ad(\d(a))x=a^ix\qua\frl a\in\kk^*\}$. Hence
$$\io(\fg_i^\d)=\{\x\in\fg^*;\Ad(\d(a))\x=a^i\x\qua\frl a\in\kk^*\}.$$
If $j\in\ZZ$, $j\ne-i$ and $\x\in\io(\fg_i^\d),x\in\fg_j^\d$ then for $a\in\kk^*$ we have
$$\x(x)=a^{-i}(\Ad(\d(a))\x)(x)=a^{-i}\x(\Ad(\d(a)\i)x)=a^{-i}a^{-j}\x(x)=a^{-i-j}\x(x)$$
hence $\x(x)=0$. Thus $\io(\fg_i^\d)\sub\fg_i^{*\d}$. Since $(\io(\fg_i^\d))$, $(\fg_i^{*\d})$ form direct sum 
decompositions of $\fg^*$ it follows that (a) holds.

Let $\fg_2^{\d!}$ be as in \cite{\LIII, 1.2}. From the definitions we see that $\io$ induces a bijection

(b) $\fg_2^{\d!}@>\si>>\fg_2^{*\d!}$.
\nl
Let $\D\in D_G$. Using (a) we see that for any $j\in\NN$ we have 

(c) $\io(\fg^\D_{\ge j})=\fg^{*\D}_{\ge j}$.
\nl
Now $\io$ induces an isomorphism $\fg_{\ge2}^\D/\fg_{\ge3}^\D@>\si>>\fg_{\ge2}^{*\D}/\fg_{\ge3}^{*\D}$. This
induces (using (b) and the definitions) a bijection $\Si^\D@>\si>>\Si^{*\D}$ (with $\Si^\D$ as in 
\cite{\LIII, 2.3}) and a bijection $\s^\D@>\si>>\s^{*\D}$ (with $\s^\D$ as in \cite{\LIII, A.1}).

We can find $\d_0\in\fD_G,\D_0\in D_G$ such that $\d_0\in\D_0$ and $G^{\D_0}_{\ge0}$ is a Borel subgroup of $G$. 
An element $\x\in\fg^*$ is nilpotent if and only if for some $g\in G$ we have $\Ad(g)\x\in\Ann(\fg^{\D_0}_{\ge0})$
(which equals $\fg^{*\D_0}_{\ge 1}=\io(\fg^{\D_0}_{\ge1})$). We see that $\io$ restricts to a bijection 
$\cn_\fg@>\si>>\cn_{\fg^*}$. Thus, if $\Ps_\fg$ (see \cite{\LIII, A.1}) is a bijection, then $\Ps_{\fg^*}$ is a 
bijection.

\subhead 2.4\endsubhead
In this subsection we assume that $G,\fg$ are as in 1.1 or as in 1.2. In both cases we
can find an isomorphism $\fg@>\si>>\fg^*$ compatible with the $G$-actions (see 1.1(a), 1.2(b)). Since $\Ps_\fg$ is
a bijection (see \cite{\LIII, A.2}) we see from 2.3 that Theorem 2.2 holds for $G$.

\subhead 2.5\endsubhead
Let $V,(,),G,\fg,\fQ(V)$ be as in 1.3. Let $\fQ(V)@>>>\fg$, $Q\m A_Q$ be as in 1.3. 

We now fix a $\ZZ$-grading $V=\op_{i\in\ZZ}V_i$ which is $s$-good (as in \cite{\LIII, 1.4}) that is, 
$\dim V_i=\dim V_{-i}\ge\dim V_{-i-2}$ for any $i\ge0$, $\dim V_i$ is even for any even $i$ and $(V_i,V_j)=0$ 
whenever $i+j\ne0$. 
For any $i$ we set 
$V_{\ge i}=\op_{i';i'\ge i}V_{i'}$.

Let $\fQ(V)_2$ be the vector space consisting of all $Q\in\fQ(V)$ such that $A_Q(V_i)\sub V_{i+2}$ for any $i$ and
$Q|_{V_i}=0$ for any $i\ne-1$. Let $\fQ(V)_2^0$ be the set of all $Q\in\fQ(V)_2$ such that 

(i) for any even $n\ge0$, $A_Q^n:V_{-n}@>>>V_n$ is an isomorphism;

(ii) for any odd $n\ge1$, $A_Q^{(n-1)/2}:V_{-n}@>>>V_{-1}$ is injective and the restriction of $Q$ to
$A_Q^{(n-1)/2}(V_{-n})$ is a nondegenerate quadratic form.
\nl
Let $\fQ(V)_{\ge2}$ be the set of all $Q\in\fQ(V)$ such that $A_Q(V_{\ge i})\sub V_{\ge i+2}$ for any $i$ and
$Q|_{V_{\ge0}}=0$. Note that $\fQ(V)_2\sub\fQ(V)_{\ge2}\sub\fQ(V)_{nil}$. 

We show:

(a) {\it if $p\ne2$, then $\fQ(V)_2^0$ is equal to the set $S$ consisting of all $Q\in\fQ(V)_2$ such that 
$A_Q^n:V_{-n}@>>>V_n$ is an isomorphism for any $n\ge0$.}
\nl
Assume first that $Q\in\fQ(V)_2^0$. Let $n\ge0$ and let
$x\in V_{-n}$ be such that $A_Q^n(x)=0$. If $n$ is even then $x=0$ by (i). If $n$ is odd then let 
$x'=A_Q^{(n-1)/2}x\in V_{-1}$. For any $x_1\in V_{-n}$ we have $(A_Qx',A_Q^{(n-1)/2}x_1)=\pm(A_Q^nx,x_1)=0$ so 
that $x'$ is in the radical of $Q|_{A_Q^{(n-1)/2}(V_{-n})}$; this radical is $0$ using (ii) and the condition 
$p\ne2$. Thus $x'=0$ that is $A_Q^{(n-1)/2}x=0$. Using the injectivity in (ii) we see that $x=0$. Thus 
$A_Q^n:V_{-n}@>>>V_n$ is injective for any $n\ge0$ hence it is an isomorphism since $\dim V_n=\dim V_{-n}$. We see
that $Q\in S$. Conversely assume that $Q\in S$. Assume that $n\ge1$ is odd. Clearly
$A_Q^{(n-1)/2}:V_{-n}@>>>V_{-1}$ is injective. If $x\in V_{-n}$ and $A_Q^{(n-1)/2}x$ is in the radical of 
$Q|_{A_Q^{(n-1)/2}(V_{-n})}$ then $0=(A_QA_Q^{(n-1)/2}x,A_Q^{(n-1)/2}x_1)=\pm(A_Q^nx,x_1)$ for any $x_1\in V_{-n}$
hence $A_Q^nx=0$ so that $x=0$ (since $Q\in S$). Thus $Q\in\fQ(V)_2^0$.

\subhead 2.6\endsubhead
As in \cite{\LIII, 2.6}, let $\fF_s(V)$ be the set of all filtrations $V_*=(V_{\ge a})_{a\in\ZZ}$ such that

(i) $\{x\in V;(x,V_{\ge a})=0\}=V_{\ge1-a}$ for any $a$;

(ii) the obvious grading of the associated vector space $gr(V_*)=\op_aV_{\ge a}/V_{\ge a+1}$ is an $s$-good 
grading with respect to the symplectic form $(,)_0$ on $gr(V_*)$ induced by $(,)$.
\nl
(Here condition (ii) can be replaced by the condition that there exists an $s$-good grading $(V_i)$ of $V$ such 
that $V_{\ge a}=\op_{a;a\ge i}V_i$ for any $a$.)

For $V_*=(V_{\ge a})\in\fF_s(V)$ let $\z(V_*)$ be the set of all $Q\in\fQ(V)$ such that
$A_Q(V_{\ge a})\sub V_{\ge a+2}$ for any $a\in\ZZ$, $Q|_{V_{\ge0}}=0$ and such that the element
$\bQ\in\fQ(gr(V_*))_2$ induced by $Q$ (see below) satisfies $\bQ\in\fQ(gr(V_*))_2^0$. The element $\bQ$ is 
defined as follows. Let $x=\sum_ax_a$ where $x_a\in V_{\ge a}/V_{\ge a+1}$. Let $\dx_a\in V_{\ge a}$ be a 
representative of $x_a$. Then $\bQ(x)=Q(\dx_{-1})+\sum_{a\le-2}(A_Q\dx_a,\dx_{-a-2})$. Note that 
$\z(V_*)\sub\fQ(V)_{nil}$.

Assuming that $p\ne2$ we note that $Q\m A_Q$ defines a bijection 

(a) $\z(V_*)@>\si>>\ti\x'(V_*)$ 
\nl
(with $\ti\x'(V_*)$ as in \cite{\LIII, A.3}); the inverse map associates to $A\in\ti\x'(V_*)$ the quadratic form
$Q:V@>>>\kk$ given by $Q(x)=(Ax,x)/2$.

We have the following result.

\proclaim{Proposition 2.7} The map $\sqc_{V_*\in\fF_s(V)}\z(V_*)@>>>\fQ(V)_{nil}$, $Q\m Q$ is a bijection.
\endproclaim
When $p\ne2$ this follows from \cite{\LIII, A.3(a)} using the bijection 2.6(a) and the bijection $Q\m A_Q$ of 
$\fQ(V)_{nil}$ onto $\cn_{\fg}$. The proof for $p=2$ will be given in 2.10, 2.11.

\subhead 2.8\endsubhead
In this subsection we assume that $V\ne0$ and that $p=2$.  
For any $Q\in\fQ(V)_{nil}$ let $e=e_Q$ be the smallest integer $\ge1$ such that $A_Q^e=0$ and let $f=f_Q$ be the 
smallest integer $\ge0$ such that $Q(A_Q^f(x))=0$ for all $x\in V$. Define a subset $H_Q\sub V$ as follows.

$H_Q=\{x\in V;A_Q^{e-1}x=0\}$ if $e\ge2f+1$;

$H_Q=\{x\in V;A_Q^{e-1}x=0,Q(A_Q^{f-1}(x))=0\}$ if $e=2f$;

$H_Q=\{x\in V;Q(A_Q^{f-1}(x))=0\}$ if $e\le2f-1$.
\nl
We show: 

$e\le 2f+1$.
\nl
It is enough to show that $A_Q^{2f+1}x=0$ for all $x\in V$ or that $(A_Q^{2f+1}x,y)=0$ for all $x,y\in V$ or that
$(A_Q^{f+1}x,A_Q^fy)=0$ for all $x,y$ or that $Q(A_Q^f(x+y))-Q(A_Q^fx)-Q(A_Q^fy)=0$ for all $x,y$; this is clear.

Let $V=\op_iV_i$ be an $s$-good grading of $V$. Let $m\ge0$ be the largest integer such that $V_m\ne0$. Let 
$Q\in\fQ(V)_{\ge2}$. We set $e=e_Q,f=f_Q$, $A=A_Q$. Let $\bQ\in\fQ(V)_2^0$. We set $\be=e_{\bQ},\baf=f_{\bQ}$,
$\bA=A_{\bQ}$. We asume that $C:=A-\bA$ satisfies $C(V_{\ge i})\sub V_{\ge i+3})$ for all $i$ and that 
$(Q-\bQ)|_{V_{\ge-1}}=0$. We show:

(i) {\it If $m$ is even then $\be=m+1$, $2\baf\le m$ (hence $2\baf<\be$). Since $\be\le2\baf+1$ we deduce that 
$\be=2\baf+1$ hence $2\baf=m$. If $m$ is odd then  $2\baf=m+1$, $\be\le m+1$ (hence $2\baf\ge\be$). In any case,
$H_{\bQ}=V_{\ge-m+1}$.}
\nl
We must show: 

if $m$ is even then $\{x\in V;\bA^mx=0\}=V_{\ge-m+1}$; if $m$ is odd then 
$\{x\in V;\bA^mx=0,\bQ(\bA^{(m-1)/2}x)=0\}=V_{\ge-m+1}$.
\nl
Let $x=\sum_{k\ge-m}x_k\in V$ with $x_k\in V_k$. We have $\bA^mx=\bA^mx_{-m}$ and if $m$ is odd then
$\bQ(\bA^{(m-1)/2}x)=\bQ(\bA^{(m-1)/2}x_{-m})$. Hence it is enough to show:

if $m$ is even then $x_{-m}=0$ if and only if $\bA^mx_{-m}=0$ (this is clear); if $m$ is odd then $x_{-m}=0$ if 
and only if $\bA^mx_{-m}=0$ and $\bQ(\bA^{(m-1)/2}x_{-m})=0$.
\nl
Assume that $m$ is odd. Assume that $\bA^mx_{-m}=0$ and $\bQ(\bA^{(m-1)/2}x_{-m})=0$. Let 
$y=\bA^{(m-1)/2}x_{-m}\in V_{-1}$. For any $z_{-m}\in V_{-m}$ we have 
$(y,\bA^{(m+1)/2}z_{-m})=(\bA^mx_{-m},z_{-m})=0$. Hence $y$ is in the radical of the symplectic form 
$(a,b)\m(a,\bA b)$ on $\bA^{(m-1)/2}V_{-m}$. Since $\bQ$ is nondegenerate on $\bA^{(m-1)/2}V_{-m}$ (with 
associated symplectic form $(a,\bA b)$) and $\bQ y=0$ we see that $y=0$. Since $\bA^{(m-1)/2}:V_{-m}@>>>V_{-1}$ 
is injective we deduce that $x_{-m}=0$. This proves (i).

We show:

(ii) {\it If $2n\ge m$ then $QA^n=0,Q\bA^n=0,\bQ\bA^n=0$. Hence if $m$ is even then $f\le m/2$; if $m$ is odd then
$f\le(m+1)/2$.}
\nl
To show that $QA^n=0$ it is enough to show that $Q(A^nx)=0$ whenever $x\in V_i$, $i\ge-m$ and $(A^{n+1}x,A^nx')=0$
whenever $x\in V_i$, $x'\in V_j$, $i,j\ge-m$, $i\ne j$. This follows from $A^nV\in V_{\ge0}$ and $Q_{V_{\ge0}}=0$,
$(V_{\ge1},V_{\ge0})=0$. The remaining equalities are proved in the same way.

The following statement is immediate.

(iii) {\it Let $P\in\End(V)$ be a sum of products of $n$ factors of which at least one is $C$ and remaining ones 
are $\bA$. Then $P(V_i)\sub V_{\ge i+2n+1}$ for all $i$. Hence if $n\ge m$ then $P=0$.}
\nl
We show:

(iv) {\it If $n\ge m$ then $A^n=\bA^n$.}
\nl
Indeed, $A^n=\bA^n+P$ where $P$ is as in (iii). Hence the result follows from (iii).

We show:

(v) {\it If $\be=2\baf+1$ then $A^{\be-1}=\bA^{\be-1}\ne0$, $A^{\be}=\bA^{\be}=0$; hence $e=\be$ and
$H_Q=H_{\bQ}$. If $\be=2\baf$ then $A^{\be-1}=\bA^{\be-1}\ne0$, $A^{\be}=\bA^{\be}=0$; hence $e=\be$.}
\nl
In the first case we have $\be=m+1$ hence $\be-1\ge m,\be\ge m$ and we use (iv); moreover, we have 
$e=\be=2\baf+1=m+1\ge2f+1$ (see (ii)), hence 
$$\align&H_Q=\{x\in V;A^{e-1}x=0\}=\{x\in V;A^{\be-1}x=0\}=\{x\in V;\bA^{\be-1}x=0\}\\&=H_{\bQ}.\endalign$$
In the second case we have $2\baf=m+1$ hence $\be=m+1$ and we continue as in the first case.

We show:

(vi) {\it If $n\ge0$, $2n+1\ge m$, then $QA^n=\bQ\bA^n$.}
\nl
If $2n\ge m$ then both sides are $0$, see (ii). Thus we may assume that $m=2n+1$. We have $A^n=\bA^n+P$, with $P$
as in (iii). By (iii) we have $P(V)\sub V_{\ge-m+2n+1}=V_{\ge0}$, hence $QP(V)=0$. We have 
$\bA^nV\sub V_{\ge-m+2n}=V_{\ge-1}$, $A\bA^nV\sub AV_{\ge-1}\sub V_{\ge1}$, $(V_{\ge0},V_{\ge1})=0$, hence 
$(A\bA^nV,PV)=0$. For $x\in V$ we have 
$$\align&QA^n=Q(\bA^nx+Px)=Q(\bA^nx)+Q(Px)+(A\bA^nx,Px)\\&=Q(\bA^nx)+(A\bA^nx,Px)=Q(\bA^nx).\endalign$$
Since $\bA^nx\in V_{\ge-1}$ we have $Q\bA^n(x)=\bQ\bA^n(x)$ as required.

We show:

(vii) {\it If $\be\le2\baf$ (hence $\baf>0$) then $QA^{\baf-1}=\bQ\bA^{\baf-1}\ne0$, 
$QA^{\baf}=\bQ\bA^{\baf}=0$. Hence $f=\baf$.}
\nl
In this case we have $2\baf=m+1$, $2(\baf-1)=m-1$. Hence the result follows from (vi).

We show:

(viii) {\it If $\be<2\baf$ then $e<2f$.}
\nl
In this case we have $f=\baf=(m+1)/2$. We must show that $e<m+1$. By (iv) we have $A^m=\bA^m$. We have 
$\be<2\baf=m+1$ hence $\be\le m$ and $\bA^m=0$. Hence $A^m=0$ and $e\le m$ as required.

Collecting together the results above we deduce:

(ix) {\it If $\be=2\baf+1$ then $e=\be=m+1\ge 2f+1$ (hence $e\ge 2f+1$). If $\be=2\baf$ then $e=\be=m+1$,
$\baf=f=(m+1)/2$ hence $e=2f$. If $\be<2\baf$ then $f=\baf=(m+1)/2$ and $e<2f$. In each case we have $H_Q=H_{\bQ}$
hence $H_Q=V_{\ge-m+1}$ and $m=\max(e-1,2f-1)$.}

\subhead 2.9\endsubhead
We preserve the setup of 2.8. We assume in addition that $\bA=0$. From the definitions we see that $V_i=0$ if 
$i\n\{-1,0,1\}$ and $\dim V_{-1}=\dim V_1\le1$. It follows that $A=0$. We have $m=0$ or $m=1$. If $m=0$ then 
$Q=\bQ=0$. If $m=1$ then $Q\ne0,\bQ\ne0$.

\subhead 2.10\endsubhead
We prove the injectivity of the map in 2.7 assuming that $p=2$. We argue by induction on $\dim V$. If $\dim V=0$,
the result is trivial. Assume now that $\dim V\ge1$. Let $Q\in\fQ(V)$ and let $V_*=(V_{\ge a})$, 
$\tV_*=(\tV_{\ge a})$ be two filtrations in $\fF_s(V)$ such that $Q\in\z(V_*)$ and $Q\in\z(\tV_*)$. We must show
that $V_*=\tV_*$. Let $\bQ\in\fQ(gr(V_*))_2$, $\bQ_1\in\fQ(gr(\tV_*))_2$ be the quadratic forms induced by $Q$.
Let $m\ge0$ be the largest integer such that $gr_m(V_*)\ne0$. Let $\tm\ge0$ be the largest integer such that 
$gr_{\tm}(\tV_*)\ne0$. If $\bQ=0$ then $Q=0$ hence $\bQ_1=0$; also, $V_{\ge1}=0$, $V_{\ge0}=V$, $\tV_{\ge1}=0$, 
$\tV_{\ge0}=V$; hence $V_*=\tV_*$ as desired. Thus we can assume that $\bQ\ne0$, $\bQ_1\ne0$. Hence $m\ge1$, 
$\tm\ge1$. Using 2.8(ix) we see that $V_{\ge-m+1}=H_Q=\tV_{\ge-\tm+1}$, $m=\max(e_Q-1,2f_Q-1)=\tm$. Thus, $m=\tm$ 
and $V_{\ge-m+1}=\tV_{\ge-m+1}$. We have $V_{\ge m}=\{x\in V;(x,V_{\ge-m+1})=0\}$, 
$\tV_{\ge m}=\{x\in V;(x,\tV_{\ge-m+1})=0\}$, hence $V_{\ge m}=\tV_{\ge m}$. Let
$V'=V_{\ge-m+1}/V_{\ge m}=\tV_{\ge-m+1}/\tV_{\ge m}$. Note that $V'$ has a natural nondegenerate symplectic form 
induced by $(,)$. We set $V'_{\ge a}=\text{image of $V_{\ge a}$ under }V_{\ge-m+1}@>>>V'$ (if $a\ge-m+1$), 
$V'_{\ge a}=0$ (if $a<-m+1$). We set $\tV'_{\ge a}=\text{image of $\tV_{\ge a}$ under }\tV_{\ge-m+1}@>>>V'$ (if 
$a\ge-m+1$), $\tV'_{\ge a}=0$ (if $a<-m+1$). Then $V'_*=(V'_{\ge a})$, $\tV'_*=(\tV'_{\ge a})$ are filtrations in
$\fF_s(V')$. Also $Q$ induces an element $Q'\in\fQ(V')$ and we have $Q'\in\z(V'_*)$, $Q'\in\z(\tV'_*)$. Note also
that $\dim V'<\dim V$. By the induction hypothesis we have $V'_*=\tV'_*$. It follows that $V_{\ge a}=\tV_{\ge a}$
for any $a\ge-m+1$. If $a<-m+1$ we have $V_{\ge a}=\tV_{\ge a}=V$. Hence $V_*=\tV_*$, as desired. Thus the map in
2.7 is injective. 

\subhead 2.11\endsubhead
We prove the surjectivity of the map in 2.7 assuming that $p=2$. By a standard argument we can assume that $\kk$ 
is an algebraic closure of the field $\FF_2$ with $2$ elements. We can also assume that $\dim V\ge2$. We choose an
$\FF_2$-rational structure on $V$ such that $(,)$ is defined over $\FF_2$. Then the Frobenius map relative to the
$\FF_2$-structure acts naturally and compatibly on the source and target of the map in 2.7. We denote each of 
these actions by $F$. It is enough to show that for any $n\ge1$ the map
$\a_n:(\sqc_{V_*\in\fF_s(V)}\z(V_*))^{F^n}@>>>\fQ(V)_{nil}^{F^n}$, $Q\m Q$ is a bijection. Since $\a_n$ is 
injective (see 2.10) it is enough to show that $|(\sqc_{V_*\in\fF_s(V)}\z(V_*))^{F^n}|=|\fQ(V)_{nil}^{F^n}|$. By 
1.4(b), 1.6(d), we have $|\fQ(V)_{nil}^{F^n}|=2^{n\dim V^2/2}$. It is enough to show that 

(a) $|(\sqc_{V_*\in\fF_s(V)}\z(V_*))^{F^n}|=2^{n\dim V^2/2}$. 
\nl
Now the left hand side of (a) makes sense when $\kk$ is replaced by an algebraic closure of the prime field with 
$p'$ elements where $p'$ is any prime number; when $p'\ne2$, this more general expression is equal to 
$p'{}^{n\dim V^2/2}$ since the map in 2.7 is already known to be a bijection in this case (we use also 1.4(b), 
1.6(d)). Then (a) follows from this equality by specializing $p'{}^n$ (viewed as an indeterminate) to $2^n$ 
provided that we can show that the left hand side of (a) (for general $p'$) is "universal" in the sense that it is
a polynomial in $p'{}^n$ with rational coefficients independent of $p',n$. 

We now compute the left hand side of (a) (for general $p'$). A collection $(f_a)_{a\in\ZZ}$ of integers is said to
be admissible if $f_a=0$ for all but finitely many $a$, $f_a$ is even for any even $a$, $f_a=f_{-a}$ for all $a$,
$f_0\ge f_{-2}\ge f_{-4}\ge\do$, $f_{-1}\ge f_{-3}\ge f_{-5}\ge\do$ and $\sum_af_a=\dim V$. For $(f_a)$ as above 
let $\cy_{(f_a)}$ be the set of all $V_*\in\fF_s(V)$ such that $\dim(gr_a(V_*))=f_a$ for all $a$. The left hand 
side of (a) is $\sum_{(f_a)}|\cy_{(f_a)}^{F^n}||\z(V_*)^{F^n}|$ where $V_*$ is any fixed element in 
$\cy_{(f_a)}^{F^n}$. Since each $|\cy_{(f_a)}^{F^n}|$ is "universal" it is enough to show that if 
$V_*\in\cy_{(f_a)}^{F^n}$ then $|\z(V_*)^{F^n}|$ is "universal". By a standard argument similar to the one in 
\cite{\LII, 1.5(d)} we see that $|\z(V_*)^{F^n}|=p'{}^{nd}|(\fQ(V)_2^0)^{F^n}|$ where $d$ is a quadratic 
expression in the $f_a$ (with integral coefficients independent of $n$) and $|\fQ(V)_2^0|$ is defined with respect
to a fixed $s$-good grading $(V_i)$ of $V$ which satisfies $\dim V_a=f_a$ and $F(V_a)=V_a$ for all $a$. It is 
enough to show that $|(\fQ(V)_2^0)^{F^n}|$ is "universal".
Let $s'$ be the number of pairs $(q,(U_{-1}\supset U_{-3}\supset U_{-5}\supset\do))$
where $q$ is a quadratic form $V_{-1}@>>>\kk$ defined over $\FF_{p'{}^n}$ and $U_i$ are subspaces of $V_{-1}$ 
defined over $\FF_{p'{}^n}$ such  that $\dim U_i=f_i$ and such that $q|_{U_i}$ is nodegenerate (for $i\le-1$ odd).
Let $s''$ be the number of all $(U_0\supset U_{-2}\supset U_{-4}\supset\do)$ where $U_i$ are subspaces of $V_0$ 
defined over $\FF_{p'{}^n}$ such  that $\dim U_i=f_i$ and such that $(,)|_{U_i}$ is nodegenerate (for $i\le 0$ 
even). From the definitions we have $|(\fQ(V)_2^0)^{F^n}|=s's''s_1$ where $s_1$ is "universal". It is easy to see
that $s''$ is "universal". It remains to show that $s'$ is "universal". If $f_0=0$ we have $s'=1$ and there is 
nothing to prove. We now assume that $f_0>0$. If $f_0$ is odd then $s'=t\nu(f_{-1},f_{-3},f_{-5},\do)$ where
$t$ is the number of nondegenerate quadratic forms on $V_{-1}$ and $\nu(f_{-1},f_{-3},f_{-5},\do)$ (as in
\cite{\LII, 1.2}) is "universal" by \cite{\LII, 1.2(a)}); moreover, $t$ is clearly "universal" hence $s'$ is 
"universal". If $f_0$ is even $\ge2$ then
$s'=t_1\nu^1(f_{-1},f_{-3},f_{-5},\do)+t_{-1}\nu^{-1}(f_{-1},f_{-3},f_{-5},\do)$ where
$t_1$ (resp. $t_{-1}$) is the number of nondegenerate quadratic forms on $V_{-1}$ which are split (resp.
nonsplit) over $\FF_{p'{}^n}$ and $\nu^\e(f_{-1},f_{-3},f_{-5},\do)$ (as in
\cite{\LII, 1.2}) is "universal" by \cite{\LII, 1.2(a)}); moreover, $t_1,t_{-1}$ are clearly "universal" hence 
$s'$ is "universal". Thus $s'$ is "universal" in any case. This completes the proof of the surjectivity of the map
in 2.7. 

\subhead 2.12\endsubhead
To give the $s$-good grading $V=\op_iV_i$ (see 2.5) is the same as to give an element $\d$ of $\fD_G$. We can then
identify $\fQ(V)_2$ with $\fg_2^{*\d}$ under the restriction of the bijection 1.3(a). We have:

(a) $\fQ(V)_2-\fQ(V)_2^0\sub\fg_2^{*\d}-\fg_2^{*\d!}$;

(b) $\fQ(V)_2^0\sub\fg_2^{*\d!}$.
\nl
Let $Q\in\fQ(V)_2-\fQ(V)_2^0$. To prove (a), we  must show that there exists $B\in Sp(V)$ such that $B$ fixes $Q$
and $B$ does not fix $V_{\ge i}$ for some $i$. 

Generally $x_k$ will denote an element of $V_k$. We set $A=A_Q$.

Assume first that $A:V_{-i}@>>>V_{-i+2}$ is not injective for some $i\ge2$. Then 
$A:V_{i-2}@>>>V_i$ is not surjective and since $\dim V_{i-2}\ge\dim V_i$, we see that $A:V_{i-2}@>>>V_i$ is not 
injective. We can find $e_{-i}\in V_{-i}-\{0\}$ such that $Ae_{-i}=0$. We can find $e_{i-2}\in V_{i-2}-\{0\}$ such
that $Ae_{i-2}=0$. Define $B\in\End(V)$ by 
$$B(\sum_kx_k)=\sum_{k\ne-i,i-2}x_k+(x_{-i}+(e_{i-2},x_{-i+2})e_{-i})+(x_{i-2}+(e_{-i},x_i)e_{i-2}).$$
By a computation exactly as that in the first two cases in \cite{\LIII, 1.7} we have
$$(B(\sum_kx_k),B(\sum_kx'_k))-(\sum_kx_k,\sum_kx'_k)=0.$$
Thus $B\in Sp(V)$.
We have
$$\align&QB(\sum_kx_k)-Q(\sum_kx_k)=Q(x_{-1})+\sum_{k\le-2;k\ne-i}(Ax_k,x_{-k-2})\\&+
(Ax_{-i}+(e_{i-2},x_{-i+2})Ae_{-i},x_{i-2}+(e_{-i},x_i)e_{i-2})-Q(x_{-1})\\&-\sum_{k\le-2}(Ax_k,x_{-k-2})
=-(Ax_{-i},x_{i-2})+(Ax_{-i},x_{i-2}+(e_{-i},x_i)e_{i-2})\\&=(Ax_{-i},e_{i-2})(e_{-i},x_i)=
-(x_{-i},Ae_{i-2})(e_{-i},x_i)=0.\endalign$$
Thus $B\i$ stabilizes $Q$.

We now assume that $A:V_{-i}@>>>V_{-i+2}$ is injective for all $i\ge2$ and that for some even $n\ge0$,
$A^n:V_{-n}@>>>V_n$ is not an isomorphism. Note that $n\ge2$.
As in \cite{\LIII, 1.7} we can find $e_{-n},f_{-n}$ linearly independent in $V_{-n}$ such that $A^ne_{-n}=0$, 
$A^nf_{-n}=0$. For $j\ge0$ we set $e_{2j-n}=A^je_{-n}$, $f_{2j-n}=A^jf_{-n}$. We have $e_n=0,f_n=0$. Also 
$e_m,f_m$ are linearly independent in $V_m$ if $m\le0$ is even. As in {\it loc.cit.}, for $j\in[0,n]$ we have 
$(e_{2j-n},e_{n-2j})=0$, $(f_{2j-n},f_{n-2j})=0$, $(e_{2j-n},f_{n-2j})=0$, $(f_{2j-n},e_{n-2j})=0$. Define 
$B\in\End(V)$ by 
$$\align&B(\sum_kx_k)=\sum_{k\n\{2h-n;h\in[0,n-1]\}}x_k\\&+\sum_{j\in[0,n-1]}
(x_{2j-n}+(-1)^j(f_{n-2j-2},x_{2j-n+2})e_{2j-n}\\&-(-1)^j(e_{n-2j-2},x_{2j-n+2})f_{2j-n}).\endalign$$
By a computation exactly as that in the third case in \cite{\LIII, 1.7} we have
$$(B(\sum_kx_k),B(\sum_kx'_k))-(\sum_kx_k,\sum_kx'_k)=0.$$ 
Thus $B\in Sp(V)$.
We have 
$$\align&QB(\sum_kx_k)-Q(\sum_kx_k)=Q(x_{-1})+\sum_{k\le-2;k\n-n,-n+2,\do,-2}(Ax_k,x_{-k-2})\\&
+\sum_{j=0}^{(n-2)/2}(Ax_{-n+2j}+(-1)^j(f_{n-2j-2},x_{-n+2j+2})e_{-n+2j+2}\\&
+(-1)^{j+1}(e_{n-2j-2},x_{-n+2j+2})f_{-n+2j+2},
\\&x_{n-2j-2}+(-1)^{n-j-1}(f_{-n+2j},x_{n-2j})e_{n-2j-2}+(-1)^{n-j}(e_{-n+2j},x_{n-2j})f_{n-2j-2})\\&
-Q(x_{-1})-\sum_{k\le-2}(Ax_k,x_{-k-2})\\&=
\sum_{j=0}^{(n-2)/2}(Ax_{-n+2j},(-1)^{n-j-1}(f_{-n+2j},x_{n-2j})e_{n-2j-2}\\&+
(-1)^{n-j}(e_{-n+2j},x_{n-2j})f_{n-2j-2})+\\&
\sum_{j=0}^{(n-2)/2}((-1)^j(f_{n-2j-2},x_{-n+2j+2})e_{-n+2j+2}\\&
+(-1)^{j+1}(e_{n-2j-2},x_{-n+2j+2})f_{-n+2j+2},x_{n-2j-2})\\&=
-\sum_{j=0}^{(n-2)/2}(x_{-n+2j},(-1)^{n-j-1}(f_{-n+2j},x_{n-2j})e_{n-2j}+
(-1)^{n-j}(e_{-n+2j},x_{n-2j})f_{n-2j})\\&
+\sum_{j=0}^{(n-2)/2}((-1)^j(f_{n-2j-2},x_{-n+2j+2})(e_{-n+2j+2},x_{n-2j-2})\\&
+\sum_{j=0}^{(n-2)/2}(-1)^{j+1}(e_{n-2j-2},x_{-n+2j+2})(f_{-n+2j+2},x_{n-2j-2})\\&=
\sum_{j=0}^{(n-2)/2}(-1)^j(f_{-n+2j},x_{n-2j})(x_{-n+2j},e_{n-2j})\\&+
\sum_{j=0}^{(n-2)/2}(-1)^{j+1}(e_{-n+2j},x_{n-2j})(x_{-n+2j},f_{n-2j})\\&
+\sum_{j=1}^{n/2}(-1)^{j-1}(e_{n-2j},x_{-n+2j})(f_{-n+2j},x_{n-2j})\\&
+\sum_{j=1}^{n/2}((-1)^j(f_{n-2j},x_{-n+2j})(e_{-n+2j},x_{n-2j})=
(f_{-n},x_n)(x_{-n},e_n)-(e_{-n},x_n)(x_{-n},f_n)\\&
+(-1)^{(n/2)-1}(e_0,x_0)(f_0,x_0)+(-1)^{(n/2)}(f_0,x_0)(e_0,x_0)=0\endalign$$
since $e_n=f_n=0$.
Thus $B\i$ stabilizes $Q$.

Now assume that $A:V_{-i}@>>>V_{-i+2}$ is injective for all $i\ge2$ and that for some odd $n\ge1$ 
the restriction of $Q$ to $A^{(n-1)/2}(V_{-n})$ is a degenerate quadratic form. Then we can find
$\x\in A^n(V_{-2n-1})-\{0\}$ such that $(A\x,A^n(V_{-2n-1}))=0$, $Q(x)=0$.  We can write
$\x=A^ne_{-2n-1}$ for a unique $e_{-2n-1}\in V_{-2n-1}-\{0\}$. For any $j\ge0$ we set
$e_{-2n-1+2j}=A^je_{-2n-1}\in V_{-2n-1+2j}$. Thus $e_{-1}=\x$ and $(Ae_{-1},A^n(V_{-2n-1}))=0$,
$Q(e_{-1})=0$. We show that $e_{2n+1}=0$. Indeed, 
$$(V_{-2n-1},e_{2n+1})=(V_{-2n-1},A^{n+1}e_{-1})=\pm(A^nV_{-2n-1},Ae_{-1})=0.$$
For $j\in[0,2n+1]$ we have
$$(e_{-2n-1+2j},e_{2n+1-2j})=0.$$
It is enough to note that
$$\align&(A^je_{-2n-1},A^{2n+1-j}e_{-2n-1})=\pm(A^{2n+1}e_{-2n-1},e_{-2n-1})\\&=(e_{2n+1},e_{-2n-1})=0.\endalign$$
Define $B\in\End(V)$ by 
$$\align& B(\sum_kx_k)=\sum_{k\ne-2n-1,-2n+1,\do,2n-1}x_k\\&
+\sum_{j=0}^{2n}(x_{-2n-1+2j}+(-1)^j(e_{2n-1-2j},x_{-2n+2j+1})e_{-2n-1+2j})\\&
=\sum_{k\ne-2n-1,-2n+1,\do,2n+1}x_k\\&
+\sum_{j=0}^{2n+1}(x_{-2n-1+2j}+(-1)^j(e_{2n-1-2j},x_{-2n+2j+1})e_{-2n-1+2j}).\endalign$$
(In the last line, $e_{2n-1-2j}$ is not defined when $j=2n+1$; we define it to be $0$.)
We have
$$\align&(B(\sum_kx_k),B(\sum_kx'_k))-(\sum_kx_k,\sum_kx'_k)
=\sum_k(x_k,x'_{-k})\\&
+\sum_{j=0}^{2n+1}(-1)^j(e_{2n-1-2j},x_{-2n+2j+1})(e_{-2n-1+2j},x'_{2n+1-2j})\\&
+\sum_{j=0}^{2n+1}(-1)^j(e_{2n-1-2j},x'_{-2n+2j+1})(x_{2n+1-2j}),e_{-2n-1+2j})\\&
-\sum_k(x_k,x'_{-k})
=\sum_{j=0}^{2n}(-1)^j(e_{2n-1-2j},x_{-2n+2j+1})(e_{-2n-1+2j},x'_{2n+1-2j})\\&
+\sum_{j=0}^{2n}(-1)^j(e_{-2n-1+2j},x'_{2n-2j+1})(x_{-2n+1+2j}),e_{2n-1-2j})=0.\endalign$$
Thus $B\in Sp(V)$.
We have     
$$\align&QB(\sum_kx_k)-Q(\sum_kx_k)\\&
Q(x_{-1}+(-1)^n(e_{-1},x_1)e_{-1})+\sum_{k\le-2;k\ne-2n-1,-2n+1,\do,2n-1}(Ax_k,x_{-k-2})\\&+
\sum_{j=0}^{n-1}(Ax_{-2n-1+2j}+(-1)^j(e_{2n-1-2j},x_{-2n+2j+1})e_{-2n+1+2j},\\&
x_{2n-1-2j}+(-1)^j(e_{-2n-1+2j},x_{2n-2j+1})e_{2n-1-2j})
\\&-Q(x_{-1})-\sum_{k\le-2}(Ax_k,x_{-k-2})\\&
=(Ax_{-1},(-1)^n(e_{-1},x_1)e_{-1})\\&
+\sum_{j=0}^{n-1}(-1)^j(Ax_{-2n-1+2j},e_{2n-1-2j})(e_{-2n-1+2j},x_{2n-2j+1})\\&
+\sum_{j=0}^{n-1}(-1)^j(e_{2n-1-2j},x_{-2n+2j+1})(e_{-2n+1+2j},x_{2n-1-2j})\\&
=-(-1)^n(e_{-1},x_1)(x_{-1},e_1)\\&
+\sum_{j=0}^{n-1}(-1)^{j+1}(x_{-2n-1+2j},e_{2n+1-2j})(e_{-2n-1+2j},x_{2n-2j+1})\\&
+\sum_{j=1}^n(-1)^{j-1}(e_{2n+1-2j},x_{-2n+2j-1})(e_{-2n-1+2j},x_{2n+1-2j})\\&
=-(-1)^n(e_{-1},x_1)(x_{-1},e_1)-(x_{-2n-1},e_{2n+1})(e_{-2n-1},x_{2n+1})\\&+(-1)^{n-1}(e_1,x_{-1})(e_{-1},x_1)=0
\endalign$$
since $e_{2n+1}=0$. Thus $B\i$ stabilizes $Q$.

This completes the proof of (a).

Now let $Q\in\fQ(V)_2^0$ and let $B\in Sp(V)$ be such that $Q(B(x))=Q(x)$ for all $x\in V$. To prove (b) it is 
enough to show that 

$B\in G^\d_{\ge0}$. 
\nl
The proof is similar to that in the last paragraph of \cite{\LIII, 1.8}.
We argue by induction on $\dim V$. Recall that $V\ne0$. Let $A=A_Q$. We
have $AB=BA$. Let $m$ be the largest integer $\ge0$ such that $V_m\ne0$. If $m=0$ we have $G^\d_{\ge0}=G$ and the
result is clear. Assume now that $m\ge1$. If $m$ is even we have $A^mV=V_m$, $\ker(A^m:V@>>>V)=V_{\ge-M+1}$. Since
$BA=AB$, the image and kernel of $A^m$ are $B$-stable. Hence $B(V_m)=V_m$ and $B(V_{\ge-M+1})=V_{\ge-M+1}$. Hence
$B$ induces an automorphism $B'\in Sp(V')$ where $V'=V_{\ge-m+1}/V_m$, a vector space with a nondegenerate 
symplectic
form induced by $(,)$. We have canonically $V'=V_{-m+1}\op V_{-m+2}\op\do\op V_{m+1}$ and $\fQ(V')_2$, 
$\fQ(V')_2^0$ are defined in terms of this ($s$-good) grading. Now $Q$ induces an element $Q'\in\fQ(V')_2^0$ and 
we have $Q'(B'(x'))=Q'(x')$ for all $x'\in V'$. If $V'=0$ then clearly $B\in G^\d_{\ge0}$ and the result is clear.
Hence we can assume that $V'\ne0$ and that the result holds for $V'$. We see that for any $i\in[-m+1,m-1]$, the
subspace $V_i+V_{i+1}+\do+V_{m-1}$ of $V'$ is $B'$-stable. 
Hence the subspace $V_{\ge i}$ of $V$ is $B$-stable. We see that $B\in G^\d_{\ge0}$, as required. Next we assume
that $m$ is odd. We have $V_{\ge-m+1}=\{x\in V;A^mx=0,Q(A^{(m-1)/2}x)=0\}$ (we use that $Q\in\fQ(V)_2^0$). Since
$B$ commutes with $A$ and preserves $Q$ we see that $B$ preserves the subspace
$\{x\in V;A^mx=0,Q(A^{(m-1)/2}x)=0\}$ hence $BV_{\ge-m+1}=V_{\ge-m+1}$.
We have $V_m=\{x\in V;(x,V_{-m+1})=0\}$. Since $B$ preserves the subspace $V_{\ge-m+1}$ and $B$ preserves $Q$ and
$(,)$ we see that $BV_m=V_m$. Hence $B$ induces an automorphism $B'\in Sp(V')$ where $V'=V_{\ge-m+1}/V_m$, a 
vector space with a nondegenerate symplectic form induced by $(,)$. We have canonically 
$V'=V_{-m+1}\op V_{-m+2}\op\do\op V_{m+1}$ and $\fQ(V')_2$, $\fQ(V')_2^0$ are defined in terms of this ($s$-good)
grading. Now $Q$ induces an element $Q'\in\fQ(V')_2^0$ and 
we have $Q'(B'(x'))=Q'(x')$ for all $x'\in V'$. If $V'=0$ then clearly $B\in G^\d_{\ge0}$ and the result is clear.
Hence we can assume that $V'\ne0$ and that the result holds for $V'$. We see that for any $i\in[-m+1,m-1]$, the
subspace $V_i+V_{i+1}+\do+V_{m-1}$ of $V'$ is $B'$-stable. Hence the subspace $V_{\ge i}$ of $V$ is $B$-stable. We
see that $B\in G^\d_{\ge0}$. This completes the proof of (b).

From (a),(b) we deduce

(c) $\fQ(V)_2^0=\fg_2^{*\d!}$.

\subhead 2.13\endsubhead
Assume that $G,\fg$ are as in 1.3. In this case Theorem 2.2 follows from 2.7, in view of 2.12(c).

\subhead 2.14\endsubhead
Assume that $G,\fg$ are as in 2.2. As in \cite{\LIII} let $\fU_G$ be the set of $G$-orbits on $D_G$. From 
\cite{\LIII, 2.1(b)} we see that $\fU_G=\fU_{G_\CC}$. In particular, $\fU_G$ is a finite set which depends only on
the type of $G$, not on $\kk$. For any $\co\in\fU_G$ we set
$$\cn_{\fg^*}^\co=\Ps_{\fg^*}\sqc_{\D\in\co}\s^{*\D}.$$
The subsets $\cn_{\fg^*}^\co$ are called the {\it pieces} of $\cn_{\fg^*}$. They form a partition of $\cn_{\fg^*}$
into smooth locally closed subvarieties (which are unions of $G$-orbits) indexed by $\fU_G=\fU_{G_\CC}$. Now 
assume that $\kk$ is an algebraic closure of the finite prime field $\FF_p$ and that we are given a split
$\FF_p$-rational structure on $G$. Then $\fg,\fg^*,\cn_{\fg^*}$ have induced $\FF_p$-structures and each $\co$ as
above is defined over $\FF_p$. Also, each piece $\cn_{\fg^*}^\co$ is defined over $\FF_p$ (with Frobenius map 
$F$). Let $n\ge1$. We have the following result.

(a) {\it $|\cn_{\fg^*}^\co)^{F^n}|$ is a polynomial in $p^n$ with integer coefficients independent of $p,n$.}
\nl
For type $A,D$ this follows from \cite{\LII, 1.8}. For type $C$ this follows from the proof in 2.11.

\enddocument